\input amstex
\input xy
\xyoption{all}
\documentstyle{amsppt}
\document
\magnification=1200
\NoBlackBoxes
\nologo
\hoffset1.5cm
\voffset2cm
\vsize15.5cm

\def\K{\bold{K}}
\def\M{\overline{M}}
\def\C{\bold{C}}
\def\P{\bold{P}}
\def\Q{\bold{Q}}
\def\S{\bold{S}}
\def\Z{\bold{Z}}

\def\cC{\Cal{C}}
\def\cO{\Cal{O}}

\bigskip



\bigskip

\centerline{\bf DESSINS FOR MODULAR OPERAD}

\medskip

\centerline{\bf AND  GROTHENDIECK--TEICHM\"ULLER GROUP}

\medskip

\centerline{\bf No\'emie~C.~Combe,\quad  Yuri~I.~Manin,\quad  Matilde~Marcolli}

\bigskip

\quad \quad \quad \quad \quad \quad \quad \quad \quad \quad  {\it ... j'avais retenu la fl\`eche, une cible qui ne sera jamais}

\quad \quad \quad \quad \quad \quad \quad \quad \quad \quad \quad {\it atteinte, une division \`a l'infini, le myst\`ere du diable.}

\smallskip

\quad \quad \quad  \quad \quad \quad \quad \quad \quad \quad \quad \quad \quad \quad \quad {\it Bernard Cun\'eo.``Le Chat du Typographe''.}

\medskip

{\it ABSTRACT.} A part of Grothendieck's program for studying the Galois group $G_{\Q}$
of the field of all algebraic numbers $\overline{\Q}$ emerged from his insight that
one should lift its action upon $\overline{\Q}$  to the action of $G_{\Q}$ upon the (appropriately defined)
profinite completion
of $\pi_1(\P^1 \setminus \{0,1, \infty\})$. The latter admits a good combinatorial encoding
via finite graphs ``dessins d'enfant''.

This part was actively developing during the last decades, starting with foundational
works of A. Belyi, V. Drinfeld and Y. Ihara. 

Our brief note concerns another part of Grothendieck program, in which
its geometric environment is extended to moduli spaces of algebraic curves,
more specifically,  stable curves of genus zero with marked/labelled points. 
Our main goal is to show that  dual graphs of such curves may play the
role of ``modular dessins'' in an appropriate operadic context.

\bigskip

{\it CONTENTS}

\medskip

0. Introduction and summary
\smallskip

1. Graphs and operads 

\smallskip

2. Dessins for modular operad: geometry and combinatorics

\smallskip

3. Dessins for  modular operad: Galois symmetries

\smallskip

4.  From operads to quantum statistical mechanical systems

\newpage

\centerline{{\bf 0. Introduction and summary}} 
\bigskip

An approach to the description of the (profinite completion of) ``absolute Galois group'' $G_{\Q}$
of the field of all algebraic numbers starts with observation that for any algebraic manifold 
(or more generally, integral scheme $X$)
 $G_{\Q}$ acts by outer automorphisms upon \'etale fundamental group
of  $X\otimes_{\Q}\overline{\Q}$ via exact sequence
$$
1\to \pi_1(X\otimes_{\Q}\overline{\Q}) \to \pi_1(X) \to G_{\Q} \to 1
$$
(see [Gr63], [SchGr64], and [Fr17], Ch. 12, for further details and references).

\smallskip

In the most studied case, that of $X = \P^1 \setminus \{0,1,\infty\}$,
the action of $G_{\Q}$ is further reduced to its action upon the so called
{\it dessins d'enfant}, that are finite graphs of very special origin
and structure. Each dessin is the inverse image of $[0,1]$
upon a Riemannian surface $Y$ which is finite  covering $Y\to \P^1$ ramified only over $\{0,1,\infty\}$.

\smallskip

In this article our aim consists in demonstrating that if we replace above 
$X = \P^1 \setminus \{0,1,\infty\}$ by the family of moduli spaces of
stable genus zero algebraic curves with labelled points,
the role of dessins d'enfant can be played by dual combinatorial
graphs of such curves.

\smallskip

This demonstration is the theme of the central Sections 3 and 4 of the paper,
whereas Sections 1 and 2 are preparatory. These sections introduce
a not quite standard combinatorial approach to the operadic formalism,
sufficient for our goals, but avoiding appeal to the machinery
of Quillen's model categories, that was used in the article [BriHoRo19]
and the monograph [Fr17] for similar purposes.

\smallskip

In Section 1 we describe, following [BoMa07], theory of graphs
in the categorical language, and provide details of operad theory
based on it.

\smallskip

In Section 2, we describe in more details the genus zero modular
operad in this environment and show why related combinatorial graphs
deserve their name of ``modular dessins''
as geometric objects.

\smallskip

In Section 3 we introduce a modification/enrichment of modular operad
necessary for defining the action of the absolute Galois
group upon modular dessins.

\smallskip

Finally, Section 4 transplants the powerful machinery of quantum statistics
from [BosCo95], [CoKr00] and [CoMar08] to the operadic context, and more specifically,
the context of enriched genus zero modular operad.

\smallskip

It shows, in particular, how one can extract a rich arithmetic information
about the absolute Galois group, using this machinery.

\bigskip

\centerline{\bf 1. Graphs and operads }

\medskip

All the following constructions are made in a fixed small universe.
The basic objects we will be considering in this preparatory section 
are {\it trees} and their categories, {\it operads} as functors on a category
of trees, and {\it symmetries} of these categories and functors.
\smallskip

This environment can be considered as a version of  ``dendroidal''
constructions described in  [MoeWe07],  [CiMoe13], but for the restricted purposes of
this paper we prefer to use basic definitions and constructions from [BoMa07].

\smallskip

Below we offer a very condensed survey of them. An interested reader
will find much more details in [BoMa07].

\medskip

1.1. {\bf Finite graphs, trees, and stable trees.}  In [BoMa07], Definition 1.1.1,
a combinatorial finite graph $\tau$ is defined as a family of finite sets $(F_{\tau}, V_{\tau})$ ({\it flags} and  {\it vertices})
and maps  $( \partial_{\tau} : F_{\tau} \to V_{\tau},\ j_{\tau}: F_{\tau} \to F_{\tau})$ ({\it boundary maps} and 
{\it structure involution}
satisfying $j_{\tau}^2=id$.)

\smallskip

Two--element orbits of $j_{\tau}$ form the set $E_{\tau}$ of {\it edges} of $\tau$.
Elements of one such orbit are sometimes called ``halves" of the respective edge,
and two points, -- boundaries of a member of this orbit, -- the boundary of the respective edge
itself.

\smallskip

One--element orbits of $j_{\tau}$ are called {\it tails}, or {\it  leaves} (we will use both words as
synonymous). A graph $\tau$ with one vertex and no edges is called {\it corolla}.

\smallskip

The {\it multiplicity} of a vertex $v$ is  the number of flags whose boundary is $v$.

\smallskip

Let us call a sequence of vertices $v_1, v_2, \dots ,v_n$, $n\ge 2$ of a graph $\tau$
{\it a path} connecting $v_1$ and $v_n$, if for each $i$,  $v_i \ne v_{i+1}$ constitute the boundary
of an edge.
 We say that the graph $\tau$
is {\it connected}, if any two its vertices can be connected by a path.
Corollas are also connected.

\smallskip

Graphs $\tau_1$ and $\tau_2$ are called {\it disjoint}, if $V_{\tau_1}\cap V_{\tau_2} =\emptyset$ and
$E_{\tau_1}\cap E_{\tau_2} =\emptyset$. We can define the disjoint union $\sqcup$
of any finite family of pairwise disjoint graphs in an obvious way. Clearly, each graph
is a disjoint union of its maximal connected subgraphs, that can be called
its {\it connected components.}

\smallskip

A (connected) graph is called a (connected) tree, if no vertex is connected by two different edges,
and more generally, there is no ``cyclic'' path
in it of length $\ge 2$. A tree is called {\it stable}, if each vertex is a boundary of
at least three different flags. If we say simply {\it stable tree}, without mentioning
connectedness, this means that each coonected component of this graph is a connected stable tree.
Later in Sec. 4, we will call  non--necessarily connected trees also {\it forests}.

\smallskip

In many papers involving graphs,  authors prefer to bypass a set--theoretic step,
to start directly with categorical definitions, and illustrate their constructions
by pictures.

\smallskip

To the contrary, here we are inclined to stay as long as possible
with set--theoretic (``combinatorial'')  notions. 
A passage to ``pictures'' is also described in this way  in [BoMa07], 1.1.2,
under the name  {\it geometric realisation of a graph.} 

\smallskip

Our  next goal consists in defining morphisms between graphs,
in such a way that we get a category $Graph$ whose objects are (some) graphs,
and each morphism $h: \tau \to \sigma$
is a triple $(h^F, h_V, j_h)$ of the following structure:
$h^F: F_{\sigma} \to F_{\tau}$ is a contravariant map,
$h_V: V_{\sigma} \to V_{\tau}$ is a covariant map,
and $j_{\tau}$ is an involution on the set of flags
of $\tau$ contained in $F_{\tau} \setminus h^F(F_{\sigma})$.

\smallskip

These data must satisfy a pretty long list of conditions/restrictions,
for which we refer the reader to [BoMa07], Definition 1.1.2.
Quite important is the end of this Definition,
saying that {\it composition of morphisms corresponds to
the set--theoretic compositions of $h^F$ and $h_V$},
and in addition explaining the behavior of $j$. 

\smallskip

{\it Example: isomorphisms of  graphs.} According to this
definition, any isomorphism $\sigma \to \tau$  induces a bijection $p_V:\,V_{\sigma}\to V_{\tau}$
of vertices, and a bijection  $p^F:\,F_{\tau}\to F_{\sigma}$
These two
bijections must be compatible in the following sense: if $f\in F_{\sigma}$, then
$\partial_{\tau}( (p^F)^{-1}f)  =p_V^{-1}( \partial_{\sigma}(f)).$
Finally, the induced map upon edges also must be a bijection.

\medskip

1.2. {\bf Operadic composition of  graphs.} Let $(\Cal{C}, \otimes )$
be a monoidal category. Generally, an operad in $\Cal{C}$ 
based upon a category of graphs $Graph$ is a functor
$(Graph , \sqcup ) \to  (\Cal{C}, \otimes )$ endowed with additional
structures. These structures would admit a natural description, {\it if}
the operation $\sqcup$ were a monoidal structure, or even in a more
enriched environment, if we first define upon $Graph$ a structure
of (simplicial) model category: see [Ho17], [BriHoRo19] and references
therein. 
\smallskip

However, $\sqcup$ is {\it not} a monoidal product even in the category of sets:
see a brief discussion in [BoMa07], sec. 1.6 and 1.7.

\smallskip

The working version of operads that we will adopt here, starts with
definition of what we mean by {\it operadic compositions} in $Graph$
and its subcaregory of stable trees,
and proceeds by extending it to operadic composition  in some algebraic manifolds.

\medskip

1.3. {\bf Definition.} {\it   Let $(\tau_i, t_i)$, $i=1,2$, $t_i\in F_{\tau_i}$, $t_1\ne t_2$,
be two pairs (graph, tail). Its composition (``grafting'') is the graph, denoted
$$
(\tau_1, t_1) * (\tau_2, t_2) ,
$$
or else $\tau_1 *_{(t_1,t_2)} \tau_2$.
It is the set theoretic union of $\tau_1$ and $\tau_2$, in which $t_1$
and $t_2$ are now halves of a new edge.}
\smallskip

Notice that  if $\tau_1$, $\tau_2$ are disjoint  (stable) trees, then their composition
is a (stable) tree as well.

\medskip

Below we will use only the following very restricted definition
of operad of stable trees $Tree$:

\medskip

1.4. {\bf Definition.} {\it An operad $Tree$ consists of a family of stable trees,
together with a family of binary grafting operations, and their iterations in all possible ways.}

\medskip

The reader will easily see that these iterations satisfy an appropriate
reformulation of usual operadic axioms.

\medskip

1.5. {\bf Example: magma operad $Mag$.} We will give two equivalent formalisms
describing this operad (see [BriHoMo19], Definitions 6.9 and 6.10).

\smallskip

{\it Description 1.} For a finite set $S$ of cardinality $n\ge 2$, define a finite set of stable trees
$Tree^S$, in the following way: $\tau \in Tree^S$ if and only if it is connected,
each vertex is the boundary of precisely three flags, and moreover,
there exists a linear ordering $s:=(s_1, s_2, \dots ,s_n)$ of $S$ such that:

\smallskip

(a) $F_{\tau} = \{(s_i,i, +)\,  |\,  i=1, \dots , n\} \cup \{(s_i,i, -)\, |\, i=1, \dots , n\} \cup (*_S,-)$.
Intuitively, $+$, resp. $-$ describes orientation of of flag towards, resp. outwards,
the vertex of the flag. 

\smallskip
(b) For each vertex $v\in V_{\tau}$, of three flags whose boundary is $v$, two flags
are oriented towards  $v$, and one outwards.

\smallskip

(c) Halves of each edge have opposite orientation.

\smallskip

(d) For each vertex $v$, there exists exactly one oriented path of edges
connecting $v$ to the root vertex $*_S$.

\smallskip

Finally, binary operadic compositions between $Tree^S$ and $Tree^T$ are allowed only if they produce
another three with the same properties, i.e. root of $\tau_1 \in Tree^S$ must be grafted to a
non--root of $\tau_2 \in Tree^T$.

\smallskip

We may add to this operad a degenerate tree, corresponding to $n=1$. We omit its description.

{\it Description 2.} For  a non--empty finite set $S$ of cardinality $\ge n$ define inductively
the set of nonassociative (and noncommutative) monomials $Mag_n^S$:

\smallskip

      (a) $Mag_1^S =S$.
      
\smallskip

        (b) $Mag_m^S = \cup_{p+q=m; p,q\ge 1} M_p^S\times M_q^S$ for all $2\le m \le n$.    
        
\smallskip

Elements of $Mag_n^S$ are written in [BriHoMo19], Definition 6.10, as linear words 
in the alphabet   $S\cup \{( , )\}$, and their identification with binary trees becomes
intuitively clear.

\medskip

Finally, we will define general operads in this context as follows.

\smallskip

1.6. {\bf Definition.} {\it  Let $(\Cal{C}, \otimes )$ be a monoidal category.
Then a stable tree operad  $A$ in $(\Cal{C}, \otimes )$ encoded by a tree operad $Tree$
consists of a family of objects of $\Cal{C}$ labelled by stable trees from $Tree$, and family of
binary operators between them, labelled by graftings $*_{(t_1,t_2)}$}.

\medskip

1.7. {\bf Comments.}  Let us explain, how the most standard definition of
the operad fits into our one.

\smallskip

Usually by an operad in  $(\Cal{C}, \otimes )$ one means a collection
of objects $A(n)$ of $\Cal{C}$ and a family of morphisms (operadic compositions)
$$
A(n)\otimes A(k_1)\otimes ... \otimes A(k_n) \to A(k_1+\dots +k_n -n)
$$
satisfying standard associativity conditions and some non--universal
additional data and restrictions such as structural actions of $S_n$ upon  $A(n)$,
inequalities upon $n,k_i$ etc. Usually one distinguishes them,
calling the respective modified operads cyclic ones, PROP's etc.

\smallskip

Our approach, developed in [BoMa07], insists on encoding all  initial data
and axioms for them in the definition of an appropriate category
of graphs (in our context trees) and and their morphisms. The magma operad
is a good illustration of this principle.

\smallskip

So from our viewpoint, an ordinary description of operad above means that
$A(n)$ is its value at the corolla with  $n$ tails oriented towards the vertex $v$
and one tail oriented outside (``root''). The operadic composition
mentioned above comes from grafting of $n$ such corollas to all non--tails
of the corolla with $n$ incoming tails.

\bigskip

\centerline{\bf 2. Dessins for modular operad:}

\centerline{\bf geometry and combinatorics}

\medskip

2.1. {\bf Combinatorics. } In this section we shall  introduce the central example of 
stable tree operad in the category of algebraic manifolds with direct product: {modular genus zero operad},
described here in the context of Definition 1.6.

\smallskip

Modular dessins ({\it ``dessins d'un vieillard'' })  that we will be considering here are  {\it combinatorial graphs} encoding
stable curves of genus zero with a finite subset of marked/labelled nonsingular points.  We will start working
over the field of algebraic numbers $\overline{\Q}$.
\smallskip

One such stable curve $C$  has only double points as singularities,
and all its irreducible components are isomorphic to $\P^1$.
Each its irreducible component must have $\ge 3$ points each of which
is either labelled, or singular.

\smallskip

2.2. {\bf Encoding stable curves.}   It is known
(see below) that up to deformation every such stable curve is encoded 
by the following combinatorial tree $\tau=\tau_C$:

\smallskip

$E_{\tau}$ := the set of double points of $C$.

$F_{\tau}$ := union of the set of labelled points (together with their labels) of $C$ and
halves of the edges from $E_{\tau}$. It is convenient to identify
these halves with preimages of double points on the normalization
of $C$.

$j_{\tau}$ sends each labelled point to itself, and each half of the edge
to another half.

$V_{\tau}$ := the set of irreducible components of $C$.

$\partial_{\tau} (f) = v$, if either $f$ encodes a labelled point, and $v$
encodes the irreducible component, to which this point belongs;
or else $f$ encodes a half of the edge, and $v$ encodes the
the respective irreducible component of the normalisation.

\medskip

2.3. {\bf Moduli spaces and their strata.} We will start with some
basic facts about moduli spaces of genus zero stable curves with marked points.
Our principal sources here are [Ke92], [Ka93], and their extension and generalisation 
in [BrMe13].

\smallskip

The main facts can be concisely stated as follows.

\medskip

 (a) Let $S$ be a finite set of cardinality $n+1$, $n\ge 3$.
Then stable genus zero  curves with $n+1$ points labelled by $S$
are parametrised by points of the smooth projective irreducible manifold $\overline{M}_{0,S}$
of dimension $n-2$.
\smallskip
The subspace of points corresponding to only {\it irreducible} curves is
an open Zariski dense submanifold  $M_{0,S} \subset \overline{M}_{0,S}$.
The graph of any such curve is corolla with $S$ tails.

\smallskip

(b) More generally, given a stable connected tree $\tau$ with the set of tails (labelled by) $S$, all
 stable genus zero modular curves with graph $\tau$ and their further specialisations/degenerations
are parametrised by the Zariski closed smooth projective manifold $\overline{M}_{0,\tau} \subset \overline{M}_{0,S}$.
\smallskip
Those curves whose graph is exactly $\tau$ are parametrised by the Zariski
open dense submanifold $M_{0,\tau} \subset \overline{M}_{0,\tau}$.

\medskip

We will call the submanifolds $ \overline{M}_{0,\tau}$, resp.  $M_{0,\tau}$,  {\it closed},
resp. {\it open} strata of structural stratification of  $\overline{M}_{0,S}$.

\medskip

2.4. {\bf Operadic compositions.} We can now sketch the definition
of the stable tree operad $\overline{M}$ in the monoidal category of algebraic
manifolds with direct product, in the framework of Definition 1.6 above.

\smallskip

Its objects labelled by stable trees are $\overline{M}_{0,\tau}$, and the (binary)
operadic composition is defined by the simple--minded formula
$$
\overline{M}_{0,\tau_1}*\overline{M}_{0,\tau_2} := \overline{M}_{0,\tau_1*\tau_2} ,
$$ 
where for brevity we omitted notation for tails.

\medskip

2.5. {\bf Examples.} (a) We start with strata of codimension one. 

\smallskip

Closed strata correspond to labelled trees having one edge. 
Up to isomorphism, they are classified by unordered
2--partitions $S = S^{\prime} \sqcup S^{\prime\prime}$, both parts
of each have cardinalities $\ge 2.$
Each part labels tails at one vertex of the edge.

\smallskip

(b) More generally, closed stratum $\M_{0,\sigma}$ is a substratum of another
one  $\M_{0,\tau}$ of relative codimension one, iff $\sigma$ can be obtained from
$\tau$ by inserting one extra edge in place of a vertex $v$ of $\tau$
and distributing half edges  (or tails) at $v$ according to
a two--partition as above.

\smallskip

By induction, we see that embeddings $\M_{0,\sigma}\subset \M_{0,\tau}$
of relative codimension $d\ge 1$ are classified by subsets of edges of $\sigma$
of cardinality $d$ such that their ``blowing down" produces $\tau$. In particular,
they can be obtained by iterating embeddings of codimension one. 

\smallskip

(c) Consider now dessins of strata 
having maximal codimension $n-2 = \roman{dim}\  \M_{0,S}$.
From the description in (b) one sees that
if one forgets the labelling of tails (half--edges) of such a graph, all such dessins
have the same structure: $n+1$ vertices are linearly ordered, say as
$\{v_1, \dots , v_{n+1}\}$; consecutive pairs $(v_1,v_2), (v_2,v_3), \dots , (v_q,v_{n+1})$
are connected by one edge each; finally, $v_2,\dots ,v_n$ carry one additional leaf (or tail),
whereas $v_1$ and $v_{n+1}$ carry two additional leaves each.
\smallskip

Then the labelling is simply a bijection between $n$ and the set of half--edges of $\tau$.

\smallskip

Surprisingly, our modular dessins of {\it maximal codimension}
with forgotten labelling form a subclass of dessins d'enfant that
occur also in the classical Grothendieck--Teichm\"uller context:
 they are what Grothendieck called {\it ``clean dessins''}.
 
 \smallskip
 
 In fact, in order to pass from our description to Grothendieck's one
 should do some re--encoding:
 we must introduce  extra vertices (and edges) and label the set of all
 resulting vertices
 as ``black'' ones and ``white'' ones. This operation is
 (almost) uniquely determined by the geometry of our trees:
 
 \smallskip
 
 (c1) Add one vertex at the free end of each leaf and one vertex
 in the middle of all edges $(v_i, v_{i+1}).$
 \smallskip
 (c2) Call all old vertices $v_i$ and $n+3$ new  ones {\it black vertices}.
 \smallskip
 (c3) Call $n$ new vertices in the middle of old edges {\it white ones.}
 \medskip
 
 In the initial Grothendieck's approach, the resulting ``bipartite'' graphs
 (black/white vertices) encode a subfamily
 of  Belyi maps $\Sigma \to \P^1$ with a very special ramification
 profile.
 \smallskip
 In he last subsections 4.13--4.23 of Section 4, we will show
 how modular dessins encode both geometric and Galois
 symmetries of the genus zero modular operad.

 \medskip
 
 (d) Now we can clarify somewhat the geometry of locally closed
 strata $M_{0,\tau}$.
 
 \smallskip
 
 From the definition, it follows that
 $$
 M_{0,\tau} = \M_{0,\tau} \setminus ( \bigcup_{\sigma} \M_{0,\sigma})
 $$
 where the union is taken over all substrata of relative codimension one,
 that in turn bijectively correspond to edges of $\tau$.

\medskip

2.6. {\bf Combinatorics of admissible projections.} Let now $S\subset S^{\prime}$
be a finite set and its subset. We call the respective {\it admissible projection} 
the morphism ${\M}_{0,S^{\prime}} \to {\M}_{0,S}$
forgetting points with labels in $S^{\prime}  \setminus S$.
We
consider relationships between dessins for ${\M}_{0,S^{\prime}}$ and for $ {\M}_{0, S}$.

\smallskip

(a) {\it Divisorial strata.} Divisorial strata of $ \overline{M}_{0, S}$ correspond to
(stable) 2--partitions $S=S_1\sqcup S_2$, and divisorial strata of $ \overline{M}_{0, S^{\prime}}$ 
correspond to (stable) 2--partitions $S=S^{\prime}_1\sqcup S^{\prime}_2$.

\smallskip

Under the admissible projection $p: \overline{M}_{0,S^{\prime}} \to \overline{M}_{0, S}$,
one such stratum projects to another one precisely when $p$  induces admissible projections
of components 
$$
S_i^{\prime}\sqcup \{pt_i^{\prime}\} \to S_i \sqcup \{pt_i\} ,\quad i=1,2.
$$ 
Here $pt_i, pt_i^{\prime}$ correspond to halves of  edges mentioned in 2.3 (a) above.

\smallskip

(b) {\it Strata of maximal codimension.} As in 2.3 (c) above, strata of maximal
codimension in  $ \overline{M}_{0, S^{\prime}}$ are encoded by "linear graphs"
(sequence of neighbouring vertices connected pairwise by edges)
that are stabilised by adding two labelled tails at each end of the graph
and one labelled tail at each middle vertex. The total set of labels
is $S^{\prime}$. 

\smallskip

Under an admissible projection  $p: \overline{M}_{0,S^{\prime}} \to \overline{M}_{0, S}$,
one should first delete all tails labelled by elements of $S^{\prime}\setminus S$.
After that one should contract all edges that did not occur in the respective
tree for a stratum of $\overline{M}_{0, S}$.

\medskip

2.7. {\bf The Grothendieck--Teichm\"uller monoid and group.} In the remaining part of this section,
we will describe very sketchily, following [BriHoRo19] and [Fr17], how the Grothendieck--Teichm\"uller symmetries,
first made explicit in [Dr90] and [Ih94], reappear in the context of combinatorial skeleton
of genus zero modular operad.

\smallskip

Let $\widehat{\bold{F}}_2$ be the profinite completion of the free group
with generators $x,y$, and $\widehat{\Z}$ the profinite completion of $\Z$.

\smallskip

The Grothendieck--Teichm\"uller monoid $\underline{\widehat{GT}}$
is the monoid of endomorphisms of  $\widehat{\bold{F}}_2$
of the form
$$
x\mapsto x^{\lambda}, \quad y\mapsto f^{-1}y^{\lambda}f
$$
where $(\lambda, f)\in \widehat{\Z} \times \widehat{\bold{F}}_2$
satisfy the following equations:
\smallskip
(a) $f(x,y)f(y,x)=1,$
\smallskip
(b) $ f((xy)^{-1},x)\, (xy)^{-m} f(y,(xy)^{-1})\, y^m f(x,y)\, x^m \,=\,1$, $m= (\lambda -1)/2$,
\smallskip
(c) ``pentagon relation'', whose precise form  we omit here; see our basic references.

\smallskip

The Grothendieck--Teichm\"uller group ${\widehat{GT}}$, by definition,  is the subgroup of invertible
elements of $\underline{\widehat{GT}}$.

\medskip

2.8. {\bf Embedding $G_{\Q} \to {\widehat{GT}}$.} The absolute Galois group $G_{\Q}$
acts upon dessins d'enfants, used to encode coverings $Y\to \P^1$ unramified
outside $\{0,1,\infty\}$ (cf. Section 0 above and basic references). This action
can be translated into the embedding $G_{\Q} \to {\widehat{GT}}$.

\smallskip

The problem of characterisation of the image of this embedding still remains unsolved.

\medskip

2.9. {\bf Example: braids and their encoding by graphs.} An action of 
the Grothendieck--Teichm\"uller monoid and group upon 
 combinatorial modular operad, which plays the central role
in [BriHoRo19] and [Fr17], proceeds via replacement of whole families
of objects and morphisms by their homotopical versions.
\smallskip
As a part of this replacement, several operads governing ``braiding relations''
of the type 2.7 (a), (b), (c) are defined,
in particular operads of parenthesized braids $PaB$ and
parenthesized ribbon braids $PaRB$: see [BriHoRo19],
Definitions 6.11 and 6.12.
\smallskip

In our on--going project, striving to avoid the introduction of
homotopical algebra, we have to engineer encoding 
 such braiding operations by graphs.

\smallskip

We finish this section by saying a few words about it.

\smallskip

In Example 1.5 above, we have shown how to encode orientation
of flags in a graph in order to define ``inputs'' and ``outputs''
in the combinatorial presentation of ordinary operads.

\smallskip

Here we can use the similar, albeit slightly more
sophisticated trick. In order to explain it, focus on the description
of braids as  morphisms in the Remark 6.3 of [BriHoRo19].
Cut such a braid in two halves and represent it as the composition
of two other braids: {\it input braid} and {\it output braid}.
Encode these halves by some words in a fixed finite alphabet,
and include in the notation information about
input/output, such as $+$ and $-$ in Example 1.5.

\smallskip

In the last Section 4, we will see that similar tricks are needed 
in order to introduce quantum statistical counting of
modular dessins: see 4.4.

\bigskip

\centerline{\bf 3. Dessins for modular operad: Galois symmetries}

\medskip

Although the geometric origins of Grothendieck's dessins d'enfant
and of our modular dessins are very different,  moduli spaces of curves
appeared  at a very early stage of Grothendieck's research dedicated
to Galois groups of algebraic numbers (at least 1984, or earlier): cf. [Gr97].

\smallskip

In our present context, the Galois group $Gal (\overline{\Q}/\Q)$ enters the scene
via its action upon the family of sets $\M_{0,S}(\overline{\Q})$ compatibly with
its tower structure for  appropriate ``admissible'' categories of labelling sets $S$.

\smallskip

We shall start with some preparatory considerations.

\smallskip

 Our approach here is based
upon the fact that manifolds $\M_{0,S}$ and some of the structural morphisms
between them have canonical models defined over $\Q$, and thus
also over algebraic extensions $K\supset \Q$ obtained by the scalar extensions:
see [Ka93] and [BrMe13].

\smallskip

Therefore, representations of the profinite Galois group $G$ of $\overline{\Q}/\Q$
in the automorphism groups of Kapranov models
define forms of $\M_{0,S}$ over field of algebraic numbers: see [Se13].

\smallskip

We will show that these forms can be united in an enriched genus zero modular
operad, upon which $G$ acts compatibly with operadic structure.

\medskip

3.1. {\bf Kapranov model of $\M_{0,n+1}.$} Consider a family of $n$
points $p_1, \dots , p_n$ in $\P^{n-2}$. Assume that they are in general position,
in the sense that any subfamily of $k\le n-1$ of these points spans a 
projective subspace $\P^{k-1}\subset \P^{n-2}$. 

\smallskip

Now construct the following tower of successive blow ups of $\P^{n-2}$:
first, blow up all points $p_i$; second, in the resulting manifold, blow up
all inverse images of lines, spanned by pairs of points $(p_i,p_j)$
(notice that these inverse images have empty intersections);
third, blow up inverse images of planes spanned by triples of points
$(p_i, p_j, p_k)$, and so on.

\smallskip

The upper floor of this tower will be our standard model of $\overline{M}_{0,n+1}$.
Clearly, it is defined over $\Q$, as well as the action of $\S_n$ upon it,
corresponding to all possible renumberings of $(p_1,\dots , p_n)$.
As was proved in [BrMe13], after an arbitrary extension $K$ of ground field,
the full automorphism group of $\overline{M}_{0,n+1}\otimes_{\Q} K$
remains the same.

\smallskip

Generally, for an arbitrary finite set $S$, the automorphisms of $\M_{0,S}$
act upon Kapranov models by permutation of $S$.

\medskip

3.2. {\bf $K$--forms of $\M_{0,S}$.} It follows that if $K$ is a normal
algebraic extension of $\Q$ with Galois group $G^K$, then $K$--forms of $\M_{0,S}$
are in a natural bijection with actions of $G^K$ upon $S$: see [Se13], 
Chapter III.

\smallskip

This will allow us to define a tree operad upon which there is a natural
action of the (profinite completion)  $G:=Gal (\overline{\Q}/\Q)$
and connect it with a similar extension of the genus zero modular operad.

\medskip

3.3. {\bf Stable tree operad with Galois action.}  We will now enrich our
definition of stable tree operad by declaring that 

\smallskip

{\it (a) flags of any corolla, and hence of any tree, are not just finite sets, but finite $G$--sets.

\smallskip

(b) Binary compositions (and their iterations) are allowed only if they are compatible
with actions of $G$. }

\medskip

3.4. {\bf Genus zero modular operad with Galois action.} Similarly, we
will now consider enriched modular operad, whose objects
are all possible forms of $\overline{M}_{0,S}$, and compositions
are compatible with respective actions of $G$ upon labelling sets $S$.

\smallskip

In order to better visualise 
 this definition, remark that $S$ can be easily reconstructed from 
the geometry of universal family of curves $C_{0,S}\to \M_{0,S}$.
Denote by $\sigma_i\,:\, \overline{M}_{0,S} \to C_{0,S},\, i\in S$,
the family of its structural sections. Clearly,
$$
\bigcup_{i\in S}\, \sigma_i(\overline{M}_{0,S}(K))  \subset C_{0,S}(K).
$$
Moreover, two sections with different labels $i\ne j \in S$ have empty intersection.

\medskip

3.5. {\bf Proposition.} {\it $S$ is the set of equivalence classes of those points of  $C_{0,S}(K)$
that lie on any of structure sections, modulo the equivalence relation
``lying on the same section''.}

\bigskip

3.6. {\bf Galois action upon modular dessins.} It is important
to understand in more details the Galois action that we have summarily described above.

\smallskip

Two examples are worth special consideration.

\smallskip

(a) Dessins of strata of maximal codimension: see 2.5 (c) above.

Is their Galois behaviour the same as  that Grothendieck's {\it clean dessins}
after re--encoding?

\smallskip

(b) M.~Kapranov suggested  to consider strata corresponding to tri--valent trees.
The magma operad (see 1.5 above) furnishes an obvious motivation for this:
one might conjecture that this will reproduce the embedding of $G$ into
the Grothendieck--Teichm\"uller group invoked in 2.9 above.

\medskip

Our last Section is dedicated to this task.

\bigskip

\centerline{\bf 4. From operads to  quantum statistical mechanical systems}

\medskip

In this Section, we outline a strategy that translates the operadic
setting into a quantum statistical mechanical system. We outline the main steps
of this strategy in general abstract form, but we also comment on 
subtleties and difficulties that arise when one implements this
general strategy in specific cases, such as the genus zero modular operad 
with Galois action discussed in the previous sections.  

\smallskip

Our  starting point is a stable tree operad $\cO$ in a 
monoidal category $(\cC,\otimes)$, which we can think of, as discussed in
Section~1 above, as a family of objects of $\cC$ labelled by stable trees, 
together with a family of binary operators between them labelled by the
grafting operations on trees.  We also assume that there is an action
of a group $G$ on the operad $\cO$ in the sense that $G$ acts on the
family of objects labelled by trees (through an action on the trees) 
compatibly with the grafting operations.  

\smallskip

What we want to construct out of this operad is a {\it quantum statistical
mechanical system}.  Such a system consists of a complex {\it algebra of observables}
$\Cal{A}$, represented by bounded operators on a Hilbert space $\Cal{H}$ of states,
together with a time evolution. 

\smallskip

{\it Representation by bounded operators} is denoted $\pi: \Cal{A} \to \Cal{B}(\Cal{H})$.

\smallskip

{\it A time evolution} is $1$--parameter family
of automorphisms $\sigma: \bold{R} \to Aut\, (\Cal{A})$, which is generated in 
the representation on $\Cal{H}$ by a Hamiltonian operator $H$ in the following sense:
$\pi(\sigma_t(A))=e^{itH} A\, e^{-itH}$, for all $A\in \Cal{A}$.  

\smallskip

The group $G$ of our initial data should act as {\it  symmetries}
of the quantum statistical mechanical system, namely as automorphisms of the
algebra $\Cal{A}$ that commute with the time evolution.

\smallskip
 
Further data that one derives from such a system include 
the {\it partition function} $Z(\beta)\, :=\, Tr\,(e^{-\beta H})$, with $\beta$ is the inverse
temperature, and {\it equilibrium states}: certain linear functionals on $\Cal{A}$.

\smallskip

Among equilibrium states there are {\it Gibbs states}
$$ 
\varphi_\beta(A) = \frac{1}{Z(\beta)} \, Tr(A\, e^{-\beta H}) 
$$
whenever these are defined, and more generally the {\it KMS states} at
inverse temperature $\beta$, whenever that set is non--empty. 

\smallskip
We refer the reader to [BraRob97] for the general operator
algebraic setting for quantum statistical mechanics, and to Chapter~3
of [CoMa08] for a discussion of several quantum statistical mechanical
systems of arithmetic origin. 

\smallskip
In the case where the group $G$ of
symmetries is the Galois group of an extension $\K$ of $\Q$, one would also like, as part of the quantum statistical
mechanical system construction,  to obtain an {\it arithmetic 
subalgebra} $\Cal{A}_{\Q}$ of the algebra of observables $\Cal{A}$. It should have
the following properties:
\smallskip

(a) KMS states $\varphi_\infty$ at zero temperature 
evaluated on observables in $\Cal{A}_\Q$ take values in an embedding $j(\K)$ of $\K$ in $\C$:
$$
 \varphi_\infty(\Cal{A}_{\Q})\subset j(\K). 
 $$
 
(b) KMS states $\varphi_\infty$ intertwine the action of $G$ as symmetries 
of the algebra with the action on the values in $\K$:
$$ 
\varphi_\infty \circ \gamma (A) = \gamma \circ \varphi_\infty (A), \ \ \ \forall \gamma \in G, \ \ \forall A\in \Cal{A}_{\Q} .  
$$


\smallskip

We outline in the following subsections a strategy in several steps aimed at this general construction.
We also highlight the typical technical difficulties that one expects to encounter at each step.

\medskip

4.1. {\bf Operads and commutative Hopf algebras}.  Let $\Cal{O}$ be an operad in the category $Sets$.

\smallskip
There is a general construction of
an associated commutative Hopf algebra, which we will denote here by $\Cal{A}_{\Cal{O}}$, see [ChaLiv07].
The grafting operation of trees that gives the operad structure gives rise to a coproduct on this
Hopf algebra that is closely related to the admissible cuts coproduct in the
Connes--Kreimer Hopf algebra of rooted trees, see [CoKr00], [ChaLiv07], [ChaLiv01], [LaMoer06], [Moe01].
\smallskip

Here are some details.

\smallskip

The construction of [ChaLiv07] starts with  associating some posets to an operad $\Cal{O}$ in $Sets$. 

\smallskip
Let $S$ be a finite set, and let $\Cal{T}(S)$ be the set of rooted trees with the vertex set $S$ endowed with
the operadic grafting compositions
$$ 
\Cal{T}(L) \times \prod_{\ell\in L} \Cal{T}(L_{\ell}) \to \Cal{T}(S) 
$$
where in $(t, (t_\ell))$ the trees $t_\ell$ are grafted at their root to the leaves of the tree $t$. This
operation is extended to the set $\Cal{F}(S)$ of forests with the vertex set $S$. A partial order structure
is determined by the following construction: $f \leq f^{\prime}$ in $\Cal{F}(S)$ if $f^{\prime}$ can be obtained from a subforest of $f$ by a composition map. 
Each poset obtained in this way in $\Cal{F}(S)$ has a unique minimal element and a maximal element, that is a 
rooted tree. One calls such posets ``intervals".  In particular, in the collection of posets constructed
in this way from forests in $\Cal{F}(S)$, every interval is isomorphic to a product of maximal intervals.
We denote intervals in this partial order by $[f,f^{\prime}]$. 

\smallskip

Given such a collection of posets, one can construct an associated commutative Hopf
algebra over $\Q$, the ``incidence Hopf algebra", which we denote here by $\Cal{A}_{\Cal{F}(S)}$.  It is
spanned by the isomorphism classes of products of maximal intervals. The commutative
multiplication of the Hopf algebra is the product of intervals and the coproduct is given by
$$ 
\Delta [f,f^{\prime}] =\sum_{f \leq f^{\prime\prime} \leq f^{\prime}} [f, f^{\prime\prime}]\otimes [f^{\prime\prime}, f^{\prime}].
 $$
Note that in general $\Cal{A}_{\Cal{F}(S)}$ is a free commutative algebra, but not necessarily
generated by the maximal intervals, since there can be isomorphisms of products of
maximal intervals with non--pairwise isomorphic factors and with 
different numbers of factors.  

\smallskip

It is shown in sec 6.3 of [ChaLiv07]
that the incidence Hopf algebra $\Cal{A}_{\Cal{F}(S)}$  is isomorphic to the
Connes--Kreimer Hopf algebra of rooted trees. As an algebra, this is the polynomial
algebra on the rooted trees $\tau$ in which a product of rooted trees $\tau_i$ is identified with a forest
$f=\tau_1 \sqcup \cdots \sqcup \tau_n$. The coproduct is defined on a tree using admissible cuts $C$
(and extended multiplicatively to forests):
$$
 \Delta (\tau)=\sum_{C\in {Cuts}(\tau)} \rho_C(\tau) \otimes  \pi_C(\tau), 
 \eqno(4.1)
$$
where ${Cuts}(\tau)$ is the set of admissible cuts of $\tau$. One
admissible cut $C$ is defined as a set of edges of $\tau$
that contains at most one edge
in any path from the root to a leaf (including the case of the empty set). The term
$\pi_C(\tau)$ is the forest consisting of the branches removed by the cut
and the term $\rho_C(\tau)$ is the remaining pruned rooted tree after the cut is performed.
The antipode is defined recursively by $S(1)=1$,
$$ 
S(\tau)=-m(S\otimes { id} - \iota \epsilon) \Delta(\tau) 
$$
where  $m$ is the multiplication, $\iota$ the unit and $\epsilon$ the counit. 

\smallskip

In our setting, we are considering operads $\Cal{O}$ in a monoidal category $(\Cal{C},\otimes)$.
However, for our construction of a quantum statistical mechanical system 
we still want to associate to them  Hopf algebras $\Cal{A}_{\Cal{O}}$ over $\Q$. Thus, it is
convenient to still work here with the same kind of Connes--Kreimer Hopf algebra.
Namely we let $\Cal{A}_{\Cal{O}}$ be the commutative algebra over $\Q$ generated by
the isomorphism classes $X_{\tau} := [C_{\tau}]$ of the objects $C_{\tau}$ in $\Cal{C}$ parametrised by the
trees $\tau$, with $X_f:= X_{\tau_1} \cdots X_{\tau_n}$ for a forest $f=\tau_1 \sqcup \cdots \sqcup \tau_n$.
The coproduct is modelled on (4.1):
$$
 \Delta (X_{\tau})=\sum_{C\in {Cuts}(\tau)} X_{\rho_{C(\tau)}} \otimes  X_{\pi_{C(\tau)}} .
\eqno(4.2)
$$

\medskip

4.2. {\bf Group of symmetries.} Now we want to include in the list of initial data
the action of a group $G$, and describe
 a modification of the construction of
a commutative Hopf algebra $\Cal{A}_{\Cal{O}}$ such that 
the action of  $G$ on the trees would induce an action on the Hopf algebra.

\smallskip
This a priori need not be the case in general: depending on the action, admissible cuts
for a tree $\tau$ may not map under the $G$--action to admissible cuts for 
other trees $\gamma \tau$ in the same orbit, and this can potentially create a problem with
the compatibility of the $G$--action with the coproduct (4.2), which we must
somehow avoid. Here are some details.    

\smallskip

4.3. {\bf Definition.} {\it (i) An admissible cut  of $\tau$
is called $G$--balanced cut, if for each $\gamma \in G$, 
the pair $(\gamma \rho_C(\tau), \gamma \pi_C(\tau))$
is an admissible cut of the tree $\gamma \tau$.
Denote the set of such cuts  $Cuts_{G(\tau)}$.

\smallskip

(ii)  The action $\gamma \in G$ is the
action of $\gamma$ on the tree components of the forest $\pi_C(\tau)$. The Hopf algebra
$\Cal{A}_{\Cal{O},G}$ is defined as above as a commutative algebra over $\Q$, with coproduct
$$
 \Delta (X_{\tau})=\sum_{C\in {Cuts}_{G(\tau)}} X_{\rho_C(\tau)} \otimes  X_{\pi_C(\tau)}\ .
 \eqno(4.3)
$$
}                      

\smallskip

The condition that the action of the group $G$ on the set of trees is compatible
with the grafting operations of the operad, however, gives a stronger constraint
on how $G$ transforms the trees. 

\smallskip

4.4. {\bf Lemma.} {\it The compatibility of the grafting operations of trees with the $G$-action ensures
that ${Cuts}_{G(\tau)}={Cuts}(\tau)$, hence that $\Cal{A}_{\Cal{O},G}=\Cal{A}_{\Cal{O}}$.}

\smallskip

{\it Proof.} Compatibility of the grafting operations of trees in the operad $\Cal{O}$ 
 with the $G$--action means that for all $\gamma \in G$
$$
\gamma\tau_1{ *}_{(\gamma  t_1, \gamma t_2)} \gamma\tau_2 = \gamma \cdot \tau_1{*}_{(t_1,t_2)} \tau_2 \, .
\eqno(4.4)
$$
Given an admissible cut $C\in {Cuts}(\tau)$, the tree $\tau$ is given by a grafting
$$ 
\tau = \rho_C(\tau) {*}_{(\ell_i, r_i)} \pi_C(\tau), 
$$
where ${*}_{(\ell_i, r_i)}$ means the successive grafting of the root $r_i$ of the $i$-th 
component of the forest $\pi_C(\tau)$ to the $i$-th leaf $\ell_i$ of the tree $\rho_C(\tau)$.
Under the action of $\gamma \in G$, the compatibility (4.4) of the grafting
operations implies that
$$
 \gamma \tau = \gamma \rho_C(\tau) {*}_{(\gamma \ell_i, \gamma r_i)} \gamma \pi_C(\tau). 
 $$
This shows that $(\gamma \rho_C(\tau),  \gamma \pi_C(\tau))$ is indeed an admissible cut of
$\gamma \tau$. The admissibility is guaranteed by the fact that each component in the
forest $\gamma \pi_C(\tau)$ has its root $\gamma r_i$ grafted to a leaf $\gamma \ell_i$
of the tree $\gamma \rho_C(\tau)$ rather than to a leaf of one of the 
previously attached component of $\gamma \pi_C(\tau)$. 

\smallskip

Thus, the compatibility requirement (4.4) 
generally is very strong. But one can cope with it by encoding
action of $G$ by appropriate labellings of flags. Then this action will not change
combinatorics of trees themselves.

\smallskip
4.5. {\bf A semigroup action.}
The next part of the construction of a quantum statistical mechanical system
associated to the operad $\Cal{O}$ is a semigroup $\Cal{S}$ acting by endomorphisms of
the commutative algebra $\Cal{A}_{\Cal{O}}$.  Here we will require only that
$\Cal{S}$ acts by commutative algebra homomorphisms,
$\Cal{S} \subset {Hom}_{{Alg}_{\Q}} (\Cal{A}_{\Cal{O}}, \Cal{A}_{\Cal{O}})$.
\smallskip
We will  not require a 
compatibility of this semigroup action with the coproduct of $\Cal{A}_{\Cal{O}}$. 
\smallskip

There is a natural candidate for such a semigroup given an operad $\Cal{O}$.
Indeed $\Cal{O}(1)$ of an operad is always a semigroup with multiplication
given by the operadic composition $\Cal{O}(1)\otimes \Cal{O}(1) \to \Cal{O}(1)$.
Moreover, the semigroup $\Cal{O}(1)$ also acts on the components
$\Cal{O}(n)$ of the operad, for $n\geq 2$, again by the operadic composition
$\Cal{O}(1)\otimes \Cal{O}(n) \to \Cal{O}(n)$. In terms of trees, the semigroup $\Cal{O}(1)$ 
corresponds to  linear trees with one root and one leaf and the grafting of 
the leaf of one tree to the root of the next. The action of $\Cal{O}(1)$ upon  $\Cal{O}(n)$ 
is similarly determined by grafting the leaf of a linear tree to the root of a
tree with $n$ leaves. Thus, we obtain the following. 

\smallskip

4.6. {\bf Lemma.} {\it
The action of the semigroup $\Cal{S}=\Cal{O}(1)$ on the algebra $\Cal{A}_{\Cal{O}}$ is
given on the generators by $X_\tau \mapsto X_{\ell {*} \tau}$, where
$\ell$ is a linear tree and $\ell {*} \tau$ is the grafting of a leaf of $\ell$ to
the root of the tree $\tau$. The action is extended multiplicatively to
forests. }

\smallskip

Notice that an admissible cut of a linear tree cuts just one edge.
Therefore  $\pi_{C}(\ell)$ is also a linear tree and $\pi_{C}(\ell){*} \tau$ then denotes the grafting of its
leaf to the root of $\tau$.

\smallskip
4.7. {\bf Remark.}
If there is an action of a group $G$ on the set of trees, compatible
with the grafting operations of the operad $\Cal{O}$, for a reason that will become more
transparent  below, instead of considering the semigroup
$\Cal{S}=\Cal{O}(1)$ we will only consider the subsemigroup $\Cal{S}^G$ of the invariants
with respect to the $G$--action in $\Cal{O}(1)$.  

\smallskip

4.8. {\bf Remark.} In the case of the modular operad  
with stable components $\{ \overline{M}_{0,S}\}$, an implementation of
a semigroup playing the role of $\Cal{O}(1)$ in our environment
may be defined with help of
the moduli spaces $\overline{L}_n$ of pointed curves studied in [LoMa00].

\medskip

4.9. {\bf A semigroup crossed product algebra.} Starting with
a semigroup action $\alpha: \Cal{S} \to {End}(\Cal{A})$  
on an algebra $\Cal{A}$, 
we will define  a semigroup crossed product algebra in the following way.

\smallskip
 Assume that all $\alpha_{\ell}$ are invertible
on their range, and  denote by $\beta_{\ell}$ the partial inverses given by $\beta_{\ell}(\alpha_{\ell}(A))=A$ 
and  $\beta_{\ell}(A)=0$ if $A\neq \alpha_{\ell}(A^{\prime})$ for some $A^{\prime}\in \Cal{A}$. 
Note that this condition of invertibility on the range is satisfied for the action of linear trees 
$\ell$ by grafting $\ell {*} \tau$ on trees $\tau$.

\medskip

4.10. {\bf Definition.} {\it Let $\alpha: \Cal{S} \to {End}(\Cal{A})$ be an action of a semigroup $\Cal{S}$ by endomorphisms 
of an algebra $\Cal{A}$ with partial inverses $\beta_{\ell}$.
The semigroup crossed product algebra $\Cal{A} \rtimes \Cal{S}$ 
is generated by $\Cal{A}$ and by elements $S_{\ell},   S_{\ell}^*$ 
satisfying the following relations: 

(i) $S_{\ell} S_{\ell^{\prime}}= S_{\ell^{\prime} \ell}$, for all $\ell,\ell^{\prime}\in \Cal{S}$.

(ii) $S_{\ell}^* S_{\ell}=1$, for all $\ell \in \Cal{S}$.

(iii) $S_{\ell}^* \, A \, S_\ell = \alpha_{\ell} (A)$, for all $\ell\in \Cal{S}$ and all $A\in \Cal{A}$ .

(iv) $S_{\ell}\, A \, S_{\ell}^* =\beta_{\ell}(A)$, for all $\ell\in \Cal{S}$ and all $A\in \Cal{A}$.}

\smallskip

The first two conditions mean that the $S_\ell$ 
define a representation by isometries of the opposite semigroup $\Cal{S}^{op}$.

\smallskip

In the case of $\Cal{A}_{\Cal{O}}$ with the action of the semigroup of linear trees by grafting on trees, we have 
$\alpha_{\ell} X_{\tau} = X_{\ell{*}\tau}$ . The partial inverse $\beta_{\ell}$ can be applied to an element $X_{\tau}$ when the tree $\tau$ starts at the root with a linear tree $\ell$, that is, when $\tau = {\ell}{*}\tau^{\prime}$ for some other tree 
$\tau^{\prime}$, in which case we have $\beta_{\ell}(X_{\tau})= X_{\tau^{\prime}}$. One can use the notation 
$\beta_{\ell} (X_{\tau}) = X_{\ell^{-1}{*} \tau}$,
where  $\ell^{-1} {*} \tau = \tau^{\prime}$ if $\tau = \ell{*} \tau^{\prime}$, and $\beta_{\ell} (X_{\tau}) =0 $ otherwise.

\medskip

4.11. {\bf Remark.} The action of the group $G$ on the algebra $\Cal{A}_{\Cal{O}}$ extends to an action
on the crossed product $\Cal{A}_{\Cal{O}} \rtimes \Cal{S}$, where the action on the $\Cal{S}$ 
part of the crossed product algebra is trivial since we are assuming that
the semigroup consists of the elements of $\Cal{O}(1)$ that are fixed by $G$.

\medskip

4.12. {\bf Constructing a Hilbert space representation.}
Let $\Cal{H}=\ell^2(\Cal{S})$ be the Hilbert space of square integrable functions on the semigroup $\Cal{S}$
endowed with the discrete topology. We can represent the elements $S_{\ell}$ of the crossed
product algebra as bounded operators on $\Cal{H}$ by
$$
 S_{\ell} \, \epsilon_{\ell^{\prime}}=\epsilon_{\ell^{\prime}{*} \ell}, 
 \eqno(4.5) 
$$
where $\{ \epsilon_\ell \}_{\ell \in \Cal{S}}$ is the standard orthonormal basis of $\Cal{H}$, and
$\ell^{\prime}{*} \ell$ is the multiplication in $\Cal{S}$ given by the grafting of the leaf of $\ell^{\prime}$ to the
root of $\ell$. Similarly, we let the elements $S_\ell^*$ act as
$$
 S_{\ell}^*\, \epsilon_{\ell^{\prime}} = \epsilon_{\ell^{\prime\prime}} 
 \eqno(4.6)
 $$
if $\ell^{\prime\prime}{*} \ell=\ell^{\prime}$,
and $0$ otherwise.
\smallskip
For grafting of linear trees condition that $\ell''{*}\ell=\ell'$ is satisfied whenever $\ell'$ has more
nodes than $\ell$ and it's equivalent to $\ell''$ and $\ell$ being the two parts of an admissible cut of $\ell'$.

\medskip

4.13. {\bf Lemma.} {\it 
The operators $S_{\ell}$ and $S_{\ell}^*$ acting as in (4.5) and (4.6) are
isometries, and they define a representation on $\Cal{H}$ of the opposite semigroup $\Cal{S}^{op}$.
}

 {\it Proof.}
The operators $S_\ell$ and $S_\ell^*$ acting as in (4.5) and (4.6) are clearly
bounded operators on $\Cal{H}$ satisfying the isometry condition $S_{\ell}^* S_{\ell}=1$, for all $\ell\in \Cal{S}$. 
The composition satisfies the identity $S_{\ell_1} S_{\ell_2}=S_{\ell_2 \ell_1}$, hence the operators define
a representation of the opposite semigroup $\Cal{S}^{op}$.

\smallskip

We then need to construct a representation of the algebra $\Cal{A}_{\Cal{O}}$ on $\Cal{H}=\ell^2(\Cal{S})$
that is compatible with (4.5) and (4.6) through the crossed product relation.
To this purpose, and so that the construction we make here would be suitable for the
goals that we will discuss later, we now focus more specifically on the
case where the group $G$ acting on trees and grafting operations is the Galois group
$G= Gal(\overline{\Q}/\Q)$. 

\medskip

4.14. {\bf Definition.}  {\it
Let $\Cal{O}$ be an operad with an action of $G=Gal\,(\overline{\Q}/\Q)$ on the set of trees
compatibly with the grafting operations. Let $\Cal{A}_{\Cal{O}}$ be the commutative Hopf algebra 
constructed as above. Let $Hom_{Alg_\Q}(\Cal{A}_\Cal{O}, \overline{\Q})$ be the set
of commutative algebra homomorphisms from $\Cal{A}_{\Q}$ to $\overline{\Q}$. Since
$\Cal{A}_{\Cal{O}}$ is a Hopf algebra, this set is a group 
$\Cal{G}(\overline{\Q})=Hom_{{Alg}_{\Q}}(\Cal{A}_{\Cal{O}}, \overline{\Q})$, 
where $\Cal{G}$ is the
affine group scheme dual to the Hopf algebra. An element $\varphi \in \Cal{G}(\overline{\Q})$
is called a {\rm balanced character} if it intertwines the $G$--action on $\Cal{A}_{\Cal{O}}$ and the $G$--action
on $\overline{\Q}$ in the following sense:
$$
 \varphi \circ \gamma = \gamma \circ \varphi, \ \ \  \forall \gamma \in G. 
 $$
We say that a character $\varphi\in \Cal{G}(\overline{\Q})$ is bounded if (under a fixed choice of
an embedding $\overline{\Q} \hookrightarrow \C$) it satisfies the condition $| \varphi(A) |< C$ for some
$C>0$ and for all $A\in \Cal{A}_{\Cal{O}}$. }

\smallskip
The following example should convince the reader that the set of bounded balanced 
characters of $\Cal{A}_{\Cal{O}}$ is non-empty.

\smallskip

4.15. {\bf Example.} {\it
Let $G$ act on the set of trees of $\Cal{O}$ so that all the orbits are finite sets.
Let $\{ \tau \}$ be a set of representatives of the $G$--orbits 
on trees of $\Cal{O}$, such that the corresponding orbit has size $d_\tau := card\, Orb\, (\tau )$. Choose a set $\{ \lambda_\tau \}$ 
of algebraic numbers such that $| \gamma \lambda_\tau | \leq 1$ for all $\gamma \in G$ 
(in a fixed embedding $\bar\Q \hookrightarrow \C$), with $\deg(\lambda_\tau)=d_\tau$.
Then setting $\varphi(X_{\gamma \tau})=\gamma \lambda_\tau$ defines a bounded
balanced character.}

\smallskip

In order to see it, for each representative $\tau$ with size of the corresponding $G$--orbit $d_\tau= card\, Orb\,(\tau)$,
choose first an algebraic number $\lambda_\tau$ such that $card\, Orb\, (\lambda_\tau)=[\Q(\lambda_\tau):\Q]=\deg(\lambda_\tau)=d_\tau$. It is always possible to divide such $\lambda_\tau$ by a sufficiently large
integer so that all the $G$--orbits are contained inside the unit disk, so we can assume that this property is
satisfied for $\lambda_\tau$. Then we obtain a character $\varphi\in \Cal{G}(\overline{\Q})$ which is both
bounded and balanced.


\medskip

Now we construct a representation of $\Cal{A}_\Cal{O}$ on $\Cal{H}=\ell^2(\Cal{S})$ in the following way.
\medskip

4.16. {\bf Lemma.} {\it
Let $\varphi\in \Cal{G}(\overline{\Q})$ be a bounded balanced character as in Definition 4.14. 
Then setting
$$
\pi_\varphi(X_\tau) \, \epsilon_{\ell} = \varphi(X_{\ell {*} \tau}) \, \epsilon_{\ell}
\eqno(4.7)
$$
defines a representation of the algebra $\Cal{A}_{\Cal{O}}$ by bounded operators on the Hilbert
space $\Cal{H}=\ell^2(\Cal{S})$. Together with (4.5) and (4.6), this determines
a representation of the crossed product algebra $\Cal{A}_{\Cal{O}} \rtimes \Cal{S}$ on $\Cal{H}$.}
\smallskip

{\it Proof.}  Since $\varphi$ is a homomorphism in $Hom_{Alg_{\Q}}(\Cal{A}_{\Cal{O}}, \overline{\Q})$
and $\Cal{S}$ acts by algebra endomorphisms, 
we have 
$$
\pi_{\varphi}(X_{\tau} X_{\tau^{\prime}}) \epsilon_{\ell} =
\varphi(X_{\ell {*} \tau}) \varphi(X_{\ell {*} \tau^{\prime}}) \epsilon_{\ell} 
=\pi_{\varphi}(X_{\tau})\pi_{\varphi}(X_{\tau'}) \epsilon_{\ell} \ .
$$ 
The property that
$\pi_\varphi(X_\tau)$ is a bounded operator follows from the boundedness property of
the character. We have
$$ 
S_{\ell} \, \pi_\varphi(X_\tau) \, S_{\ell}^* \, \epsilon_{\ell'} =
\varphi(X_{\ell^{\prime\prime} {*}\tau}) \epsilon_{\ell^{\prime}} 
$$
if $\ell{*}\ell^{\prime\prime}=\ell^{\prime}$, and
$0$ otherwise.
Furthermore
$$
S_{\ell}^* \, \pi_{\varphi}(X_{\tau}) \, S_{\ell}\, \epsilon_{\ell^{\prime}} = \pi_{\varphi}(X_{\ell{*} \tau}). 
$$
\smallskip

This gives $S_{\ell}\, \pi_{\varphi}(X_{\tau}) \, S_{\ell}^*=\pi_{\varphi}(\beta_{\ell}(X_{\tau}))$ and
$S_{\ell}^* \, \pi_{\varphi}(X_{\tau}) \, S_{\ell}=\pi_\varphi(\alpha_{\ell}(X_{\tau}))$. Thus, the
relations of the semigroup crossed product algebra $\Cal{A}_{\Cal{O}}\rtimes \Cal{S}$ are satisfied.
\medskip

4.17. {\bf Remark.}
Let $\Cal{A}_{\Cal{O},\C}=\Cal{A}_{\Cal{O}}\otimes_{\Q} \C$. 
The representation $\pi_{\varphi}: \Cal{A}_{\cO} \to \Cal{B}(\ell^2(\Cal{S}))$
of Lemma 4.16 extends to a representation of $\Cal{A}_{\Cal{O},\C}$ and of the crossed product
$\Cal{A}_{\Cal{O},\C}\rtimes \Cal{S}=(\Cal{A}_{\Cal{O}}\rtimes \Cal{S})\otimes_{\Q} \C$. Let $\Cal{B}_{\Cal{O}}$ denote the
$C^*$-algebra obtained from the crossed product algebra $\Cal{A}_{\Cal{O},\C}\rtimes \Cal{S}$ by including
adjoints of the elements of $\pi_{\varphi}(\Cal{A}_{\Cal{O},\C})$ and completing it in the operator norm 
of the algebra of bounded operators
$\Cal{B}(\ell^2(\Cal{S}))$. 

\smallskip
The $C^*$-algebra $\Cal{B}_{\Cal{O}}$ is the semigroup crossed product $C^*$--algebra
$\Cal{B}_{\Cal{O}}=\Cal{A}_{\Cal{O},\pi_{\varphi}}\rtimes \Cal{S}$, where $\Cal{A}_{\Cal{O},\pi_{\varphi}}$ is the $C^*$-completion
of $\pi_{\varphi}(\Cal{A}_{\Cal{O},\C})$ in $\Cal{B}(\ell^2(\Cal{S}))$.
The $\Q$--algebra $\Cal{B}_{\Cal{O}}^{ar} :=\Cal{A}_{\Cal{O}}\rtimes \Cal{S}$
is then referred to as the {\it arithmetic subalgebra} of $\Cal{B}_{\Cal{O}}$.

\medskip

4.18. {\bf Semigroup homomorphisms and time evolution.} Now we pass to the
construction of a time evolution operator on the $C^*$--algebra 
$\Cal{B}_{\Cal{O}}=\Cal{A}_{\Cal{O},\pi_{\varphi}}\rtimes \Cal{S}$ of Remark 4.17.
The central requirement is that the time evolution commutes with the
symmetries given by the action of the group $G$.
We use the following strategy.

\smallskip

4.19. {\bf Proposition.} {\it
Suppose that there exists a semigroup homomorphism $\lambda: \Cal{S} \to \bold{N}$
to the multiplicative semigroup of natural numbers $\bold{N}$, 
with the property that the growth of the multiplicities $a_n=\{ \ell \in \Cal{S}\,|\, \lambda(\ell)=n\}$
is such that the Dirichlet series   
$$
\sum_{n\geq 1} a_n \, n^{-\beta}
\eqno(4.8)
$$
converges for sufficiently large $\beta$. Then setting
$$
\sigma_t(X_{\tau})=X_{\tau} \ \ \ \and  \ \ \  \sigma_t(S_{\ell})=\lambda(\ell)^{it}\, S_{\ell}, \ \ \  \sigma_t(S_{\ell}^*)=
\lambda(\ell)^{-it}\, S_{\ell}^*, 
\eqno(4.9)
$$
defines a time evolution on the $C^*$--algebra $\Cal{B}_{\Cal{O}}=\Cal{A}_{\Cal{O},\pi_{\varphi}}\rtimes \Cal{S}$ that
commutes with the action of $G$ by automorphisms. In
the representation of Lemma 4.16 this time evolution is generated by the
Hamiltonian
$$
H\, \epsilon_{\ell} = \log \lambda(\ell) \, \epsilon_{\ell},
\eqno(4.10)
$$
with partition function
$$
Z(\beta)=Tr(e^{-\beta H}) = \sum_{\ell \in \Cal{S}} \lambda(\ell)^{-\beta}
\eqno(4.11)
$$
that converges for sufficiently large inverse temperature $\beta$.
}
\smallskip

{\it Proof.} Since $\lambda: \Cal{S} \to \bold{N}$ is a semigroup homomorphisms with
values in a commutative semigroup, we have $\sigma_t(S_{\ell}S_{\ell^{\prime}})=
\sigma_t(S_{\ell^{\prime}\ell})=\lambda(\ell)^{it}\lambda(\ell^{\prime})^{it} S_{\ell^{\prime}\ell}$.
The action (4.9) is moreover compatible with the relations of
the semigroup crossed product algebra. 

\smallskip
Since the time evolution is trivial
on the subalgebra $\Cal{A}_{\Cal{O},\C}$ and on $\Cal{A}_{\Cal{O},\pi_{\varphi}}$, and
nontrivial on the $\Cal{S}$ part of the crossed product, while 
the action of $G$ by automorphisms is nontrivial on the $\Cal{A}_{\Cal{O},\pi_{\varphi}}$
part and trivial on the semigroup $\Cal{S}$ part, the two actions on the
crossed product algebra $\Cal{B}_{\Cal{O}}=\Cal{A}_{\Cal{O},\pi_{\varphi}}\rtimes \Cal{S}$ commute.

\smallskip

To see that the time evolution is generated by the Hamiltonian (4.10) we need to
check that, for all $A\in \Cal{A}_{\Cal{O},\C}\rtimes \Cal{S}$, 
$$ 
\pi_{\varphi}(\sigma_t(A)) =e^{it H} \pi_{\varphi}(A) e^{-itH}. 
$$
This is the case for $A\in \Cal{A}_{\Cal{O},\C}$ where the time evolution acts trivially.
\smallskip
For $S_{\ell}$ we have $e^{it H} S_{\ell} e^{-it H}\epsilon_{\ell^{\prime}}=
\lambda(\ell^{\prime} {*}\ell)^{it} \lambda(\ell^{\prime})^{-it} S_{\ell} \epsilon_{\ell^{\prime}}$, hence $e^{it H} S_{\ell} e^{-it H}=\lambda(\ell)^{it} S_{\ell} =\sigma_t(S_{\ell})$, and
similarly for $S_\ell^*$. Under the assumption that the multiplicities grow at most polynomially, $m_n\leq P(n)$, the partition function is given by the Dirichlet series (4.8)
$$ 
Z(\beta) =\sum_{n\in \bold{N}} a_n\, n^{-\beta}, 
$$
and converges for sufficiently large inverse temperature $\beta$.
\medskip

4.20. {\bf Lemma.} {\it
Let $G$ be an action on the trees of the operad $\Cal{O}$, compatible with grafting, such that
the orbits of $G$ on the set of trees are finite. Let $\Cal{S}$ be the semigroup of linear trees
fixed by the $G$--action, with the composition by grafting tail to root. Then there is a choice
of a semigroup homomorphism $\lambda: \Cal{S} \to \bold{N}$ satisfying the properties of Proposition 4.19,
such that the partition function $Z(\beta)=\sum_{\ell \in \Cal{S}} \lambda(\ell)^{-\beta}$ is convergent for
all $\beta>0$.}

\smallskip

{\it Proof.}  
In the case of the semigroup of linear trees $\ell \in \Cal{S}$
fixed by the $G$--action, an example of a
homomorphism $\lambda: \Cal{S} \to \bold{N}$ satisfying the growth condition of Proposition 4.19
can be obtained in the following way. Let $L(\ell) := card\, E(\ell)$ be the length of the linear graph $\ell$ counted as
number of edges. We assign $L(\ell)=0$ to a graph $\ell$ consisting of a single vertex, which
we include as the unit of the semigroup. Let $\Cal{L}$ be the set of labels of vertices and edges on which the group $G$
also acts. Since we are assuming that the semigroup $\Cal{S}$ consists of elements that are fixed
by $G$, the labels of such linear graphs must be in the subset $\Cal{L}^G$ of $G$--fixed points in $\Cal{L}$.
Let $k= card\, \Cal{L}^G$ be the cardinality of this set, which we assume finite, since we work under the
assumptions that orbits of $G$ on the set of trees of the operad $\Cal{O}$ are finite.

\smallskip
Let $N\in \bold{N}$ be chosen so that $N > 2 k^2$. Then setting $\lambda(\ell)=N^{L(\ell)}$
defines a semigroup homomorphism $\lambda: \Cal{S} \to \bold{N}$ to the multiplicative semigroup
of positive integers, since the length is additive under grafting: $L(\ell {*} \ell^{\prime})=L(\ell)+L(\ell^{\prime})$.
\smallskip
A linear tree $\ell$ has $L(\ell)$ edges and $L(\ell)+1$ vertices. If edges and vertices are
labelled  by $\Cal{L}^G$, this gives $k^{L(\ell)} \cdot k^{L(\ell)+1}$ possible choices of $\ell$.
Thus, the partition function $Z(\beta)$ is computed by the series
$$ 
\sum_{L\in \bold{N}} k^{L} \cdot k^{L+1} N^{-L}. 
$$
In view of the choice of $N$, this is bounded by
$$ 
\sum_{L\in \bold{N}} k^{2L +1} k^{-2L} 2^{-L} =k \sum_{L\geq 1} 2^{-L} \leq k. 
$$
Thus, in this case the partition function is convergent for all $\beta>0$. 
\smallskip
Moreover, we see that in this situation all KMS states are Gibbs states of the form
$$ 
\phi_\beta(A) =\frac{1}{Z(\beta)}  Tr(\pi_{\varphi} (A) e^{-\beta H}). 
$$
We discuss the zero--temperature ground states and their properties in the next subsection.
\medskip

4.21. {\bf Gibbs states and ground states.}
The ground states at zero temperature are defined in Chapter 3 of [CoMa09]
as the weak limits of the KMS--states at large inverse temperature $\beta$, when
$\beta\to \infty$:

$$
\phi_{\infty}(A)= \lim_{\beta \to \infty} \phi_{\beta}(A)= \lim_{\beta \to \infty} \frac{1}{Z(\beta)} 
Tr(\pi_{\varphi} (A) e^{-\beta H}) , 
$$
for Gibbs states $\phi_{\beta}$ at inverse temperature $\beta$.  The group $G$ acts
on the set of Gibbs states at a given $\beta$ by pullback, $\gamma^*\phi_{\beta}(A)=\phi_{\beta}(\gamma(A))$.
\medskip

4.22. {\bf Proposition.} {\it
Consider the time evolution of Lemma 4.20. The ground states, when
restricted to the arithmetic subalgebra $\Cal{B}^{ar}_{\Cal{O}}$, take values in $\overline{\Q}$ and satisfy the
intertwining property with respect to the $G$--action
$$ 
\phi_{\infty}\circ \gamma  = \gamma \circ \phi_{\infty} . 
$$
}
\smallskip

{\it Proof.} We consider the time evolution defined by the semigroup $\lambda: \Cal{S} \to \bold{N}$
with $\lambda(\ell)=N^{L(\ell)}$ discussed in Lemma 4.20.

The kernel of the corresponding Hamiltonian $H \epsilon_{\ell} = \log \lambda(\ell) \epsilon_{\ell}
=L(\ell) \log(N)\, \epsilon_{\ell}$ is spanned by a single vector $\epsilon_{\ell}$, corresponding to the graph $\ell$ 
consisting of a single vertex, which we have included as unit of the semigroup. We write this
vector as $\epsilon_1$.  

We then see from the above that we have
$$
 \phi_{\infty}(A) = \langle \epsilon_1, A \, \epsilon_1 \rangle, 
 $$
namely the limit is the ground state in the usual sense of projection onto the kernel of the Hamiltonian.

Consider the case where $A$ is an element of the arithmetic subalgebra $\Cal{B}^{ar}_{\Cal{O}}$. It suffices
to consider the elements $A\in \Cal{A}_{\Cal{O}}$ since any element in $\Cal{B}^{ar}_{\Cal{O}}$ that is not contained
in $\Cal{A}_{\Cal{O}}$ will have projection $\langle \epsilon_1, A \, \epsilon_1 \rangle=0$. In fact, there would
be a number of $S_\ell$ or $S_\ell^*$ terms that map $\epsilon_1$ to some other $\epsilon_{\ell^{\prime}}$
that is orthogonal to $\epsilon_1$. Thus, we can restrict ourselves to the case $A=X_{\tau}$ in  $\Cal{A}_{\Cal{O}}$. The
case of a forest is similar.
 
Since in the construction we have chosen $\varphi \in \Cal{G}(\overline{\Q})$ to be a balanced character,
we obtain 
$$ 
 \langle \epsilon_1, \pi_{\varphi}(X_{\tau}) \, \epsilon_1 \rangle = \varphi(X_{\tau}) .
 $$
Hence $\phi_{\infty}$ evaluated on $\Cal{B}^{ar}_{\Cal{O}}$ takes values in $\overline{\Q}$. We then have
$$
\phi_{\infty}(\gamma X_{\tau})=\varphi(\gamma X_{\tau}) = \gamma \, \varphi(X_{\tau})=\gamma \, \phi_{\infty}(X_{\tau}), 
$$
so we obtain the intertwining condition.

\medskip

We can also write more explicitly the Gibbs states at finite values of the inverse temperature as follows.

\smallskip

4.23. {\bf Corollary.} {\it
The values of Gibbs states at inverse temperature $\beta >0$ on elements $X_{\tau}$ are given by
$$
 \phi_\beta(X_\tau) = \frac{1}{Z(\beta)} \sum_{\ell \in \Cal{S}} \varphi(X_{\ell\star\tau})\, \lambda(\ell)^{-\beta}.
 $$  }
\smallskip

{\it Proof.} We have
$$
\phi_{\beta}(X_{\tau}) =\frac{1}{Z(\beta)} Tr (\pi_{\varphi} (X_{\tau}) e^{-\beta H}) 
 =\frac{1}{Z(\beta)} \sum_{\ell\in \Cal{S}} \langle \epsilon_{\ell}, \pi_{\varphi} (X_{\tau})\, \epsilon_{\ell} \rangle \, 
\lambda(\ell)^{-\beta} 
$$
$$
= \frac{1}{Z(\beta)} \sum_{\ell\in \Cal{S}} \langle \epsilon_1, S_{\ell}^*\, \pi_\varphi (X_{\tau})\, S_{\ell} \, \epsilon_1 \rangle \, 
\lambda(\ell)^{-\beta} 
 =  \frac{1}{Z(\beta)} \sum_{\ell\in \Cal{S}} \varphi(X_{\ell{*}\tau}) \, \lambda(\ell)^{-\beta}. 
 $$
Thus, we can regard these as  normalised generating series of the values $\varphi(X_{\ell{*}\tau})$
weighted by $\lambda(\ell)^{-\beta}$.

\bigskip

{\bf Acknowledgements.}  N. C. Combe acknowledges support from  the Minerva Fast track grant from 
the Max Planck Institute for Mathematics in the Sciences, in Leipzig.

M.  Marcolli  acknowledges support
from NSF grant DMS--1707882 and
NSERC grants RGPIN--2018--04937 and
RGPAS--2018--522593.

\bigskip 
\centerline{\bf References}

\medskip

[BoMa07]  D.~Borisov, Yu. Manin. {\it Generalized operads and their inner cohomomorhisms .}
 In: Geometry and Dynamics of Groups
and Spaces (In memory of Aleksander Reznikov). Ed. by M. Kapranov et al.
Progress in Math., vol. 265 (2007), 
Birkh\"auser, Boston, pp. 247--308.
arXiv math.CT/0609748

\smallskip

[BraRob97]  O.~Bratteli, D.W.~Robinson. {\it Operator algebras and quantum statistical mechanics. 2. Equilibrium states. Models in quantum statistical mechanics}. Second edition.Texts and Monographs in Physics. Springer-Verlag (1997). xiv+519 pp.

\smallskip

[BosCo95] J.B.~Bost, A. Connes. {\it Hecke algebras, type III factors and phase transitions with spontaneous symmetry breaking in number theory.} Selecta Math. (N.S.) 1 (1995), no. 3, pp. 411--457.
\smallskip

[BriHoRo19] P. de Brito, G. Horel, M. Robertson.  {\it Operads of genus zero curves
and the Grothendieck--Teichm\"uller group.}  Geometry and Topology,
23 (2019), pp. 299 --346. 

\smallskip

[BrMe13] A. Bruno, M. Mella. {\it The automorphism group of $\overline{M}_{0,n}$}. 
Journ. Eur. Math. Soc., 15 (2013), pp. 949--968.

\smallskip

[ChaLiv01]  F.~Chapoton, M.~Livernet, {\it Pre--Lie algebras and the rooted trees operad}, 
Internat. Math. Res. Notices,  no. 8 (2001), pp. 395--408.

\smallskip

[ChaLiv07] F.~Chapoton, M.~Livernet. {\it Relating two Hopf algebras built from an operad}.
Int. Math. Res. Not. IMRN, no. 24 (2007), Art. ID rnm131, 27 pp. arXiv:0707.3725.

\smallskip

[CiMoe13] D.--Ch. Cisinski, I. Moerdijk. {\it Dendroidal Segal spaces and $\infty$--operads.}
Journ. of Topology, 6 (2013), pp. 675 --704.

\smallskip
[CM19-1] N. Combe, Yu. Manin. {\it Genus zero modular operad and absolute Galois group.}
arXiv:math.AG/1907.10313v2. 27 pp.

\smallskip

[CM19-2] N. Combe, Yu. Manin. {\it Symmetries of genus zero modular operad}.
arXiv:math.AG/1907.10317.

\smallskip

[CoKr00]  A.~Connes, D.~Kreimer. {\it Renormalization in quantum field theory and the 
Riemann--Hilbert problem I:
The Hopf algebra structure of graphs and the main theorem}. Comm. Math. Phys., Vol.~210 (2000) N.1, 249--273.

\smallskip

[CoMar08] A.~Connes, M.~Marcolli. {\it Noncommutative Geometry, Quantum Fields and Motives}.
Colloquium Publications, AMS 55 (2008).

\smallskip

[Dr90] V.  Drinfeld. {\it On quasitriangular quasi--Hopf algebras and on a group that is
closely connected with} $Gal({\overline{\Q}}/\Q).$ (Russian). Algebra i Analiz 2 (1990),
pp. 149 --181.

\smallskip

[Fr17]  B. Fresse.  {\it Homotopy of operads and Grothendieck--Teichm\"uller groups,
Part I: The algebraic theory and its topological background.}
Math. Surveys and Monographs 217, AMS (2017), 532 pp.
\smallskip

[Gr63] A. Grothendieck. {\it Rev\^{e}tements \'etales et groupe fondamental. Fasc. I:
 Expos\'es 1 \`a 5.} S\'eminaire de G\'eom\'etrie Alg\'ebrique, vol.1960/1961, Paris (1963).

\smallskip

[Gr97] A. Grothendieck. {\it Esquise d'un Programme (1984)}. In ``Geometric Galois Actions''
(L. Schneps and P.  Lochak eds.) London Math. Soc. Lect. Notes Ser. 242.
Cambridge UP, 1997, pp. 5--47.

\smallskip

[Ho17]  G. Horel. {\it Profinite completion of operads and the Grothendieck--Teichm\"uller
group.} Adv. Math. 321 (2017), pp. 326--390.

\smallskip

[Ih94] Y.  Ihara. {\it On the embedding of  $Gal (\Q/Q)$ into $\widehat{GT}$.}
In [SchGrT94], pp. 289--305.

\smallskip
[Ka93] M.~Kapranov. {\it Veronese curves and Grothendieck--Knudsen moduli spaces
$\overline{M}_{0,n}$.} Journ. Alg. Geom., 2 (1993), pp. 239--262.

\smallskip

[Ke92]  S. Keel. {\it Intersection theory of moduli spaces of stable n--pointed curves
of genus zero.} Trans. AMS, 30 (1992), pp. 545--584.

\smallskip
[KoMa94] M. Kontsevich, Yu. Manin. {\it Gromov--Witten classes, quantum cohomology,
and enumerative geometry.}  Comm. Math. Phys., 164 (1994), pp. 525--562.

\smallskip

[LaMoe06] P.~van der Laan, I.~Moerdijk. {\it 
Families of Hopf algebras of trees and pre-Lie algebras}. 
Homology Homotopy Appl. 8 (2006) N.1, pp. 243--256.

\smallskip

[LoMa00]  A.~Losev, Yu.~Manin, {\it New moduli spaces of pointed curves and pencils of flat connections}.  Michigan Math. J. 48 (2000), pp. 443--472.

\smallskip

[Ma99] Yu.~Manin. {\it Frobenius manifolds, quantum cohomology, and moduli spaces.}
AMS Colloquium Publications, Vol. 47  (1999), xiii + 303 pp.

\smallskip
[MaMar19-1]  Yu. Manin,  M.  Marcolli.  {\it Nori diagrams and persistence homology.}
Math. Comput. Sci. 14 (2020), no. 1, pp. 77--102.
arXiv:1901.1031. 

\smallskip
[MaMar19-2]  Yu. Manin,  M.  Marcolli.  {\it Quantum statistical mechanics of the
absolute Galois group.} arXiv:mathAG/1907.13545v2. 62 pp. To appear
in SIGMA.
\smallskip

[MayMo19] J.~Maya, J.~Mostovoy. {\it Simplicial equations for the moduli space
of stable rational curves.} arXiv:math/AG/1906.052I3. 7pp.

\smallskip

 [Moe01] I.~Moerdijk. {\it
On the Connes-Kreimer construction of Hopf algebras}. In ``Homotopy methods in algebraic topology (Boulder, CO, 1999),
 Contemp. Math., AMS 271 (2001), pp. 311--321.
 
\smallskip

[MoeWe07]  I. Moerdijk, I. Weiss. {\it Dendroidal sets.} Alg. Geom. and Topology,
7 (2007), pp. 1441--1470.

\smallskip

[SchGrT94] L. Schneps, ed. {\it The Grothendieck Theory of Dessins d'enfants.} London
Math. Soc. LN Series 200 (1994), Cambridge UP.

\smallskip

[Se13] J.--P.  Serre. {\it Galois cohomology.} Springer (2013).

\bigskip

{\bf No\'emie C. Combe, Max-Planck-Institut f\"ur Mathematik in den Naturwissenschaften
Inselstr. 22, 04103 Leipzig, Germany}

\medskip

{\bf Yuri I. Manin, Max-Planck-Institut f\"ur Mathematik, Vivatsgasse 7, 53111 Bonn, Germany}

\medskip

{\bf Matilde Marcolli, Math. Department, Mail Code 253-37, Caltech, 1200 E.California Blvd.,
 Pasadena, CA 91125, USA}

\enddocument